\documentclass[12pt]{amsart}
\usepackage{amscd,amssymb}

\newtheorem{theorem}{Theorem}[section]
\newtheorem{corollary}[theorem]{Corollary} 

\newtheorem{proposition}[theorem]{Proposition}
\theoremstyle{definition}

\theoremstyle{remark}

\newtheorem{example}[theorem]{Example}

\numberwithin{equation}{section}

\newcommand{\abs}[1]{\lvert#1\rvert}
\def\norm#1{\left\Vert#1\right\Vert}
\def\R {{\mathbb{R}}}

\def\C {{\mathbb{C}}}

\def\N{{\mathbb{N}}}
\def\e{{\varepsilon}}
\def\Lip{{\mathrm{Lip}\,}}
\def\Z {{\mathbb{Z}}}
\def\diam{{\mathrm{diam}\,}}
\def\supp{{\mathrm{supp}\,}}

\begin{document}

\title[Universality property]
{On a universality property of some abelian Polish groups}

\author[S. Gao and V. Pestov]
{Su Gao and Vladimir Pestov}
\address{S.G.: Department of Mathematics, PO Box 311430,
University of North Texas, Denton, TX 76203-1430, USA}
\email{sgao@unt.edu}
\urladdr{http://www.math.unt.edu/$^\sim$sgao/}

\address{V.P.: School of Mathematical and Computing Sciences,
Victoria University of Wellington, P.O. Box 600, Wellington,
New Zealand}
\email{vova@mcs.vuw.ac.nz}
\urladdr{http://www.mcs.vuw.ac.nz/$^\sim$vova}
\address{{\rm After July 1, 2002:}
Department of Mathematics and Statistics,
University of Ottawa, Ottawa, Ontario, K1N 6N5, Canada.
}



\begin{abstract} 
We show that every 
abelian Polish group is the topological factor-group of a 
closed subgroup of the full
unitary group of a separable Hilbert space with the 
strong operator topology.
It follows that all orbit equivalence relations induced by 
abelian Polish group actions
are Borel reducible to some orbit equivalence relations 
induced by actions of the
unitary group.
\end{abstract}

\maketitle 

\section{Introduction}

For a class $\mathcal C$ of topological groups, there are usually two 
competing notions of universality. In some context, a universal object 
is a topological group for which every group in the class $\mathcal C$ 
can be isomorphically embedded as a topological
subgroup. In a different sense, a universal object means a topological 
group of which every group in $\mathcal C$ is a topological factor-group, 
i.e., there is a continuous and open
homomorphism from the universal group onto each group in $\mathcal C$. 
The notions
are sometimes distinguished from each other by being respectively 
called {\it injective universality} and {\it projective universality},
but the terminology has not been standardized. 
In either of the two senses,
it is of definite interest whether a universal object belongs to 
the class $\mathcal C$, although
the mere existence of universal objects, no matter in $\mathcal C$ 
or not, can be more 
important. 

Here we consider a universality property 
that combines the above two
senses. A topological group $G$ is universal for 
$\mathcal C$ in our sense if every group in 
$\mathcal C$ is a topological factor-group of a topological 
subgroup of $G$. This is a
weaker notion than either one mentioned above. 
It is also easy to see that the relation ``$H$ is a 
a topological factor-group of a topological subgroup of $G$''
is transitive.

The class $\mathcal C$ we deal with in
this paper is that of all abelian Polish groups. 
The existence of both injectively and projectively universal
Polish groups is already known.
However, those universal groups are of a special
kind, not at all well understood, while 
the groups universal in our weaker sense count among
them such familiar objects as 
the additive group of the Banach space $\ell_1$. Since the latter
topological group embeds into the 
full unitary group $U_{\infty}$ of
the separable complex Hilbert space $\ell_2$, equipped with the strong
operator topology, it follows that $U_{\infty}$ is universal in
our sense 
for the class of all abelian Polish groups. 
%
%
%

Our investigation is motivated by questions in descriptive set theory of 
equivalence relations. Let us briefly review the main concepts of this theory.
If a
Polish group $G$ acts in a Borel manner on a standard Borel space $X$
(in which case $X$ is called a {\it Borel $G$-space}),
we denote the induced orbit equivalence relation by $E^X_G$. 
If $E$ and $F$ are equivalence relations on standard Borel 
spaces $X$ and $Y$ respectively,
then we say that $E$ is {\it Borel reducible to} $F$, denoted $E\leq_B
F$, if
there is a Borel function $f:X\to Y$ such that, for all $x,y\in X$,
$$ xEy\iff f(x)Ff(y). $$
An important open problem in the theory of equivalence relations is: Is
there
an orbit equivalence relation induced by a Polish group action which is
not
Borel reducible to any orbit equivalence relation of an action of the
unitary group?
We provide a partial answer as follows.

\vskip .3cm
\noindent{\bf Theorem \ref{main}. }
{\it Let $G$ be an abelian Polish group 
and $X$ be a Borel $G$-space.
Then there is a Borel $U_\infty$-space $Y$ such that 
$E^X_G\leq_B E^Y_{U_\infty}$.}
\vskip .3cm

The following interesting question seems to be open:
Is every separable metrizable topological group a
topological factor-group of a suitable topological subgroup
of $U_\infty$? If the answer to this question is in the affirmative, 
then the abovementioned
open problem about orbit equivalence relations would 
be completely settled.

A by-product of our investigation is 
a new (and more elegant) proof of the known result from 
\cite{M-P}:
every separable metrizable abelian topological group embeds as
a topological subgroup into a monothetic metrizable topological
group.

The two main tools used in our paper are
transportation distances and positive definite functions. 
Transportation distances have been independently
discovered in different areas of mathematics and are thus
known under numerous names. We give a survey of the theory
in Section \ref{transp}. Section \ref{pdf} outlines the use
of positive definite functions to construct strongly continuous
unitary representations of some topological groups.
In section \ref{borel} the main Borel reducibility results are deduced.
We have attempted to make the
article relatively self-contained, collecting in it definitions and
hints of proofs of known results for reader's convenience.

\section{Transportation distances\label{transp}}

Transportation distances were initially introduced by Kantorovich
in his 1942 paper \cite{Kant}
in order to study the classical mass transportation problem,
and have since then found numerous applications in different areas of
mathematics, in some of which they have been rediscovered
independently and explored to varying degrees of depth and from
various angles.

\subsection{Free normed spaces}
Let $X=(X,d,\ast)$ be a pointed metric space, that is, a triple 
where
$d$ is a metric on a set $X$ and $\ast\in X$ is a
distinguished point. Denote by $L(X,\ast)$, or simply
by $L(X)$, the real vector space having $X\setminus\{\ast\}$
as its Hamel basis and $\ast$ as zero. 


There obviously exists the
largest prenorm, $p$, on $L(X)$ with the property that
the distance induced on $X$ does not exceed $d$: for all
$x,y\in X$, $p(x-y)\leq d(x,y)$. 

Such a $p$ is in fact a norm, and
the restriction of the associated distance to $X$ coincides with $d$. 
Indeed, these are equivalent to saying that every metric
space isometrically embeds into a normed space as a 
linearly independent set. Here is such an embedding
(described in \cite{K-R} and, independently, \cite{AE}, cf.
also \cite{mich}.)
Denote by $\Lip(X,\ast)$ the Banach space of all Lipschitz
functions $f\colon X\to\R$ with the property $f(\ast)=0$,
where $\norm f$ equals
the smallest Lipschitz constant for $f$. 
For an $x\in X$, denote by $\hat x$ the evaluation
functional:
\[\Lip(X,\ast)\ni f\mapsto f(x)\in \R.\]
The mapping 
\[X\ni x\mapsto \hat x\in \Lip(X,\ast)^\prime\]
is an isometric embedding of $X$ into the dual Banach space of
$\Lip(X,\ast)$ as a linearly independent subset (an easy check). 

In fact, more is true: every element of $L(X,\ast)$, if considered as 
a finitely-supported measure on
$X\setminus\{\ast\}$, determines a bounded linear functional
on $\Lip(X,\ast)$, and thus $L(X,\ast)$ embeds into the
dual Banach space $\Lip(X,\ast)^\prime$ as a normed subspace.
The dual norm on 
$L(X)$ induced from $\Lip(X,\ast)^\prime$ is exactly the maximal
prenorm that we are after. Notice also that $X$ is closed in
$L(X)$. 

The normed space $L(X)$ has the following universal property,
which provided the main motivation for such investigations as
\cite{AE,Rai,Fl1,Fl2}.

\begin{theorem}
\label{univer}
Let $E$ be a normed space, and let $f\colon X\to E$ be a
1-Lipschitz map with the property $f(\ast)=0$. Then there is
a unique linear operator $\bar f\colon L(X)\to E$ of norm $1$
extending $f$. 
\end{theorem}

\begin{proof} The existence of a unique linear operator
$\bar f$ as above is clear. It remains to notice that
the prenorm on $L(X)$ denoted by $q(x)=\norm{\bar f(x)}_E$
has the property $q(x-y)\leq d(x,y)$ for all $x,y\in X$,
and thus $q(z)\leq\norm z$ for all $z\in L(X)$ and the statement
follows.
\end{proof}

The formula (\ref{graev}) 
that follows can be seen both as the definition
of the transportation distance \cite{Kant}, and as an alternative
description of the norm of the free normed space
\cite{AE,Rai,Fl1,Fl2} going back to Graev \cite{Gr},
where it appeared in the context of free (abelian) groups.

\begin{theorem} 
Let $x\in L(X,\ast,d)$. Then
\begin{eqnarray}
\norm x &=& \inf\left\{\sum_{i=1}^n \abs{\lambda_i}d(x_i,y_i)\colon 
n\in\N,~\lambda_i\in\R,~x_i,y_i\in X,\right.\nonumber \\
& & \left. x=\sum_{i=1}^n \lambda_ix_i,
~0= \sum_{i=1}^n \lambda_iy_i\right\}.
\label{graev}
\end{eqnarray}
\label{graev2}
\end{theorem}

\begin{proof}
Denote by $\norm\cdot^\prime$ the prenorm determined by 
the expression on the right hand side of the formula (\ref{graev}),
and let $\norm\cdot$ stand
for the norm of the free normed space $L(X)$.
If $x\in L(X)$, then for any two decompositions of $x$ and $0$
as in (\ref{graev}) one has
\begin{eqnarray*}
\norm x&\leq&\sum_{i=1}^n\norm{\lambda_i(x_i-y_i)} \\
&=& 
\sum_{i=1}^n \abs{\lambda_i}d(x_i,y_i),
\end{eqnarray*}
and consequently $\norm x\leq \norm x^\prime$.
Now let $x,y\in X$. Writing
$x-y=1\cdot x + (-1) y$ and $0 = 1\cdot x+ (-1)x$, one concludes that
\begin{eqnarray*} \norm{x-y}^\prime &\leq& 1\cdot d(x,x)+1\cdot d(x,y) 
\\
&=& d(x,y),
\end{eqnarray*}
and consequently 
$\norm x^\prime\leq\norm x$ for every $x\in L(X)$.
\end{proof}

The Banach space completion of $L(X)$ is
denoted by $B(X)$ and called the
{\it free Banach space} on the pointed metric space 
$(X,\ast,d)$. It has an universal property
similar to that in Theorem \ref{univer} with respect to
all {\it Banach} spaces $E$.

\begin{example}
If $X=\Gamma\cup\{\ast\}$ is a set equipped with a discrete 
($\{0,1\}$-valued) metric, 
the free Banach space $B(\Gamma\cup\{\ast\})$
(where $\ast$ is the distinguished point) is
isometrically isomorphic to $\ell_1(\Gamma)$.
\end{example}

It is easy to see that the following three conditions are
equivalent: (i) a metric space $X$ is separable; (ii) the free normed
space $L(X)$ is separable; (iii) the free Banach space
$B(X)$ is separable.

On this occasion let us
remind a well-known and simple fact from classical Banach space
theory. (Cf. e.g. p. 108 in \cite{LZ}.)

\begin{theorem}
Every separable Banach space $E$ is a factor-space of $\ell_1$. 
\label{ell-one}
\end{theorem}

\begin{proof}
Let $f$ be an arbitrary map from $\N_+$ onto an everywhere dense
subset of the sphere of radius $\frac 12$ around zero in $E$.
The map $f$ is 1-Lipschitz with respect to the 
$\{0,1\}$-valued metric on $\N$,
and thus extends to a linear operator $\bar f$ 
of norm $\leq 1$ (in fact, exactly $1$)
from $B(\N)\cong\ell_1$ to $E$. (Here $0\in\N$ serves as the
distinguished point.) 
Let $x\in E$ be arbitrary with $\norm x=\frac 12$.
It is possible to choose recursively a sequence of 
elements $k_n\in\N$ and non-negative 
scalars $\lambda_n\leq 2^{-n}$ in such a way
that each element $\sum_{i=1}^n\lambda_if(x_i)$ is at a distance
$<2^{n+1}$ from $x$. Consequently, $z=\sum_{i=1}^\infty\lambda_ix_i$ 
is in $\ell_1$ and $\bar f(z)=x$. Thus, the operator $\bar f$ is
onto, and the Open Mapping Theorem finishes the proof.
\end{proof}

Let $a,b\in X$. The following fact, standard in theory of
free objects, is established 
by applying the universal property from Theorem \ref{univer} to
the 1-Lipschitz mapping $X\ni x\mapsto x-a+b\in B(X,b)$.

\begin{proposition} 
Let $X=(X,d)$ be a metric space. For all choices 
of the distinguished point $\ast\in X$ the resulting
free Banach spaces
$B(X,d,\ast)$ are isometrically 
isomorphic between themselves. \qed
\label{allareiso}
\end{proposition}

Let us assume temporarily 
that $(X,d)$ has diameter $\leq 1$. Denote by
$X^\dag$ the metric space obtained from $X$ by adding an extra point
$\{\dag\}$ at a distance $1$ from every $x\in X$. 
Denote by $\phi$ the linear functional of norm $1$ on
$B(X^\dag,\dag)$ which takes
$X$ to $\{1\}$ and which exists by Theorem \ref{univer}.
Let $B(X)_0$ stand for the kernel of $\phi$, and let
$L(X)_0 = B(X)_0\cap L(X)$.

\begin{proposition} Assume that $\diam X\leq 1$.
For every $\ast\in X$, the free Banach space $B(X,d,\ast)$ 
{\rm (}respectively the free normed space $L(X,d,\ast)${\rm )}
is isometrically isomorphic to $B(X)_0$ {\rm (}respectively,
$L(X)_0${\rm )}.
\label{kernel}
\end{proposition}

Here, similarly to the proof of Prop. \ref{allareiso}, the isomorphic
embedding $B(X,d,\ast)\hookrightarrow B(X^\dag,\dag)$, when
restricted to $X$, is of the form $x\mapsto x-\ast$.

Recall that if $\mu$ is a 
measure on the product of two standard
Borel spaces $X$ and $Y$, then the {\it marginals} of $\mu$
are the push-forward measures $\pi_{i,\ast}\mu$, $i=1,2$, along
the coordinate projections. 
The (finitely-supported) signed
measures $\mu$ on $X\times X$ whose marginals are, respectively,
$x$ and $0$, can be identified with a pair of representations
of $x$ and $0$ as in (\ref{graev}).

Denote by $\tilde d$ the distance determined by the free norm,
$\tilde d(x,y) =\norm{x-y}$. 
Theorem \ref{graev2} and Proposition \ref{kernel}
lead to the following result.

\begin{theorem}
\label{marginals} 
Let $\mu_1,\mu_2$ be finitely-supported
probability measures on $X$. Then 
\begin{equation} 
\label{kantor}
\tilde d(\mu_1,\mu_2)
= \inf\left\{\int_{X\times X} d(x,y)\, d\nu~ \colon~ 
\pi_{i,\ast}(\nu) = \mu_i,~ i=1,2\right\}.
\end{equation} 
\qed
\end{theorem}

The formula (\ref{kantor}) makes sense for arbitrary Borel
probability measures on a metric space, and
the distance $\tilde d$ is known in this and similar contexts as the
{\it transportation distance,} {\it Monge--Kantorovich distance,} 
{\it Prokhorov distance} (in probability),
{\it Wasserstein distance} (in ergodic theory),
or else {\it Earth Mover's distance} (in computer science).
See the two-volume monograph \cite{RR}, largely devoted to
the study of the transportation distance and containing a very
large -- though still not exhaustive -- bibliography.

Here is the master result. (Cf. e.g. \cite{RR},
Section 4.1.)

\begin{theorem}[Kantorovich optimality criterion] 
Let $X=(X,d)$ be a metric space.
A probability measure $\nu$ on $X\times X$ whose marginals are, respectively,
$\mu_1$ and $\mu_2$, achieves the infimum in (\ref{kantor}) if and only
if there exists a 1-Lipschitz function $f\colon X\to\R$
such that for all pairs $(x,y)\in\supp\nu$ one has
\[f(x)-f(y) = d(x,y).\]
\qed
\label{criterion}
\end{theorem} 

The following is an immediate consequence.
(Cf. \cite{Rai,Fl1,Fl2} for direct proofs.) 

\begin{corollary}
The infimum in the
formula (\ref{graev}) is achieved at some representations 
of $x$ and $0$ with $x_i,y_i\in\supp x\cup\{0\}$. \qed
\label{suppo}
\end{corollary}

\begin{corollary}[Integer Value Property]
If $x$ is a linear combination 
with integer coefficients, then the infimum in 
(\ref{graev}) is achieved at some representations
of $x$ and $0$ as linear combinations with integer coefficients.
\label{integer}
\end{corollary}

\begin{proof}
The Kantorovich criterion reduces the result to the
following fact, established by an easy combinatorial argument. 
Suppose a matrix $A$ with real entries is such that 
the entries in each column and in each row add up to an integer.
Then all non-zero entries of $A$ can be replaced
with integers without altering the column-sums and row-sums of $A$.
\end{proof}  

In the language of optimization theory, Corollary \ref{integer}
says that a transportation problem with integer supply and
demand has an integer optimal solution. 
This is a classical result, to be
found in textbooks such as \cite{Sak} (Remark 10 on page 179).

For the complex free normed spaces the above
results starting with \ref{graev2}
are no longer true \cite{Fl1,Fl2}. 

\subsection{Graev metrics on free abelian groups}
The group envelope of $X$ in $L(X,\ast)$ is just the
free abelian group having $X\setminus\{\ast\}$ as the set
of free generators. We will denote it by $A(X,\ast)$ or else
simply by $A(X)$.
The restriction of the distance 
$\tilde d$, generated by the free norm, 
to $A(X,\ast)$, which we denote by $\bar d$, 
is a bi-invariant metric, and $\bar d\vert_X =d$.
An analogue of Theorem \ref{graev2} can be stated for
$\bar d$, and together with
Corollary \ref{integer} it implies the following.

\begin{corollary}
The metric $\bar d$ is the maximal among all bi-invariant metrics on
$A(X)$ whose restriction to $X$ is majorized by $d$.
\label{important}
\qed
\end{corollary}

In theory of topological groups, the metric $\bar d$ is known
as the {\it Graev metric} \cite{Gr}.

\begin{corollary} The metric group $A(X,d,\ast)$ is a {\rm (}closed\,{\rm )}
metric subgroup of the normed space $L(X,d,\ast)$. \qed
\label{closed1}
\end{corollary}

The above two Corollaries are just equivalent forms
of the same result, first stated by Tkach\-enko \cite{T}, who 
had offered a direct,
albeit incomplete, proof. Uspenskij \cite{U} later
noted that the result in question
follows from the Integer Value Property.

The metric space completion of 
the group $A(X)$ equipped with the metric $\bar d$ is an abelian
topological group, which we will denote by
$\widehat A(X,d,\ast)$. 

\begin{corollary}
The complete metric group
$\widehat A(X,d,\ast)$ is a {\rm (}closed\,{\rm )} metric subgroup of
the Banach space $B(X,d,\ast)$.
\label{closed}
\qed
\end{corollary}

The following universal property of the metric group $A(X,d,\ast)$
is a standard result in the theory.

\begin{proposition}[Graev] Let $X=(X,d)$ be a metric space, and let
$G$ be an abelian group equipped with a bi-inv\-ariant metric
$\varsigma$, and let $f\colon X\to G$ be a 1-Lipschitz map
(with regard to the distances $d$ on $X$ and $\varsigma$ on $G$),
taking $\ast$ to $0_G$. Then there exists a unique 
1-Lipschitz homomorphism
$\bar f\colon A(X,d,\ast)\to (G,\varsigma)$ extending $f$ from $X$.
\end{proposition}

\begin{proof}
For purely algebraic reasons, there is only one group
homomorphism $\bar f\colon A(X,d,\ast)\to (G,\varsigma)$ extending $f$.
Let us show that $\bar f$ is in fact 1-Lipschitz as well.
Define a pseudometric $\rho$ on $A(X)$ as follows: for $x,y\in A(X)$,
\[\rho(x,y):=\varsigma(\bar f(x),\bar f(y)).\]
This $\rho$ is a bi-invariant pseudometric, and the restriction of
$\rho$ to $X$ is majorized by $d$. We conclude by Corollary \ref{important}
that $\rho\leq \bar d$. But this is another way of saying that
$\bar f\colon (A(X),\bar d)\to (G,\varsigma)$ is 1-Lipschitz.
\end{proof}

\begin{corollary} Let $(G,\varsigma)$ be an abelian group equipped
with a complete bi-invariant metric. Let $f\colon (X,d)\to
(G,\varsigma)$
be a 1-Lipschitz map, taking $\ast$ to $0_G$.
Then $f$ extends to a unique continuous homomorphism
$\bar f\colon \bar A(X)\to G$, which is moreover 1-Lipschitz. \qed
\end{corollary}

Here is another elementary and well-known observation,
again going back to \cite{Gr}.

\begin{proposition}
Every
metrizable abelian topological group $G$ is a topological factor-group
of a metrizable group of the form $A(X,d,\ast)$.
If $G$ is completely metrizable, it is a quotient-group of a group
of the form $\widehat A(X,d,\ast)$. If $G$ is Polish, then the latter
group can be assumed Polish as well.
\label{factors}
\end{proposition}

\begin{proof}
Let $d$ be any bi-invariant metric on $G$ generating the topology.
Set $X=G$, $\ast=e_G$, and
consider the group $A(G,e_G,d)$. The identity map
from $G$ to itself extends to a unique 1-Lipschitz homomorphism
$i\colon A(X)\to G$ onto. 
Every $\e$-neighborhood of identity, $V_\e$, in $A(X)$
contains the $\e$-neighborhood formed within $X$, and therefore
the image $\bar f(V_e)$ has a non-empty interior in $G$
(as it contains $f(V_\e\cap X)$). It follows that $\bar f$ is
an open homomorphism. The remaining statements are obvious.
\end{proof}

Proposition \ref{factors} and Corollary \ref{closed1} together
imply:

\begin{corollary}
Every abelian metrizable group $G$ is isomorphic with a topological
factor-group of a closed subgroup of the additive group of a
normed space $E$. If $G$ is complete metrizable, then $E$ is
a Banach space. If $G$ is Polish, then $E$ is separable Banach.
\qed
\end{corollary}

The authors of the paper \cite{M-M}, 
where the above result appeared in print for the first time, 
ought to have mentioned that the Corollary had in fact
entered topological folklore shortly after the
publication of Tkach\-enko's influential work \cite{T}.

Invoking Theorem \ref{ell-one}, one obtains:

\begin{corollary} Every abelian Polish
group $G$ is isomorphic with a topological factor-group of a
closed topological subgroup of the additive group of $\ell_1$.
\label{ell-on}
\end{corollary}

\begin{proof} Let $\pi\colon\bar A(X)\to G$ be a
factor-homomorphism, and let $T\colon l_1\to B(X)$ be an
open linear operator onto. Then the complete preimage, $E$, of
$\pi^{-1}(G)$ under $T$ is a closed topological subgroup of
$\ell_1$. The restriction of $T$ to the complete preimage of a
closed set is a quotient map. Consequently,
the composition $\pi\circ (T\vert_E)$ is
an open homomorphism of topological groups (as a composition
of two open homomorphisms). 
\end{proof}

\section{\label{pdf}Positive definite functions and topological
subgroups of $U_{\infty}$}

Recall that a complex-valued function $f$ on a group $G$ is
{\it positive definite} if for every finite collection
$g_i,i=1,\ldots,n$ of elements of $G$ and every
complex scalars $\lambda_i$, $i=1,2,\ldots,n$, $n\in\N$,
\[\sum_{i,j=1}^n f(g_ig_j^{-1})\lambda_i\overline{\lambda}_j \geq 0.\]
It is a standard fact in representation theory that continuous
positive definite functions on a topological group $G$
are in one-one correspondence with 
strongly continuous cyclic representations of $G$ possessing
a (distinguished) cyclic vector
(c.f., e.g. \cite{N}, \S 30). 
For a separable metrizable group $G$, its embeddability
into $U_\infty$ is thus closely related to the existence of topology-generating
positive definite functions on $G$. 
The standard argument in fact gives the following finer
result.

\begin{theorem}
Let $G$ be a separable metrizable group and $1_G$ be its identity element.
Then $G$ is isomorphic to a topological subgroup of $U_\infty$ iff there is
a continuous positive definite function on $G$ separating $1_G$ and closed
subsets of $G$ not containing $1_G$.
\label{sbgp-U}
\end{theorem}

\begin{proof}
Let $f:G\to\C$ be a continuous positive definite function on $G$ which 
separates $1_G$
and closed subsets not containing $1_G$. 
Form the linear space 
$X$ of all complex-valued
functions on $G$ with finite support. For $x,y\in X$, let
$$ \langle x,y\rangle=\sum_{g,h\in G} f(h^{-1}g) x(g)\overline{y(h)}. $$
Let $N=\{x\in X : \langle x,x\rangle=0\}$. Then $N$ is a linear subspace of $X$ and
the sesquilinear form induces an inner product on $X/N$, making $X/N$ a pre-Hilbert space.
Let $H$ be the completion of $X/N$ under the induced norm metric. Then $H$ is a separable
complex Hilbert space. The standard representation of $G$ in $U(H)$ defined by
$$ T_gx(h)=x(g^{-1}h) $$
is easily checked to be a topological embedding of $G$ 
into $U(H)$ with the strong operator topology. 

Conversely assume that $G$ is a topological subgroup of $U(H)$ with the strong operator
topology, where $H$ is some separable complex Hilbert space. Note that for any 
$v\in H$, the function 
$$ f_v(g)=\langle g(v),v\rangle $$
is continuous and positive definite on $G$. Moreover the collection $\{ f_v : v\in H\}$
generates the topology on $G$. By separability of $G$ there is a countable subcollection
which already generates the topology of $G$. Denote this subcollection by ${\mathcal F}_0$.
Then the set ${\mathcal F}_1$ of all finite products of elements of ${\mathcal F}_0$ is 
again a collection of
positive definite functions on $G$, and this new set separates $1_G$ and closed subsets
of $G$ not containing $1_G$. Let ${\mathcal F}_1=\{f_n: n\in\N\}$. Without loss of generality we can assume that $f_n(1_G)\leq 1/2^n$ for each $n\in\N$. Finally define
$$ f(x)=\sum_{n\in\N} f_n(x). $$
Then $f$ is a continuous positive definite function on $G$ separating $1_G$ and
closed subsets of $G$ not containing $1_G$.
\end{proof}

In \cite{S} Shoenberg proved that, for $1\leq p\leq 2$, 
the function $e^{-\norm x^p}$ is positive
definite on $\ell_p$. Since this function obviously 
separates the identity from closed subsets
not containing the identity, one obtains the following.

\begin{proposition}
\label{megrel}
The additive group of each $\ell_p$, 
$1\leq p\leq 2$, is isomorphic to
a closed subgroup of $U_\infty$. 
In particular, $\ell_1$ is {\rm (}isomorphic to{\rm )} a closed subgroup
of $U_\infty$. 
\qed
\end{proposition}

These facts were noted by Megrelishvili in \cite{Meg}. 

Thus we have the following in view of Corollary \ref{ell-on}.
\begin{corollary}
Every abelian Polish group is isomorphic to a 
factor-group of a closed abelian subgroup of
$U_\infty$.
\label{univU}
\qed
\end{corollary}

In the remaining part of this section we consider the topological group
$L^0(X,U(1))$. Here $X$ is an arbitrary uncountable standard Borel space with a
non-atomic Borel measure. 
The group $L^0(X,U(1))$ consists of all (equivalence classes of)
measurable functions from $X$
into the circle rotation group $U(1)$, and the topology 
of $L^0(X,U(1))$ is that of convergence
in measure. The topological group $L^0(X,U(1))$ is the unitary group of
the abelian von Neumann algebra $L^\infty(X)$, equipped with 
the ultraweak topology. This can be also considered as the
strong operator topology with regard to the standard
representation of $L^\infty(X)$ by multiplication operators in
$L^2(X)$. 

The following is another well known theorem for which it is hard
to find a standard reference.
However, arguments that are sufficient to establish the theorem can
be found in many sources, e.g., \cite{CH} and \cite{Top}.
(The result as stated appears, for instance, in the paper \cite{G-N},
but it had been certainly known to experts
for a long time before that.)

\begin{theorem}
Let $G$ be a separable metrizable abelian group. 
Then $G$ is isomorphic to a topological
subgroup of $U_\infty$ iff $G$ is isomorphic to a 
topological subgroup of $L^0(X,U(1))$.
\label{L0}
\qed
\end{theorem}

\begin{proof}
The remarks above have given an embedding of $L^0(X,U(1))$ into $U_\infty$.
Now suppose $G$ is an abelian topological subgroup of $U_\infty$. 
We can extend $G$ to a
maximal abelian von Neumann algebra $W$. 
Then $W$ is isomorphic to $L^\infty(X)$, and the
ultraweak topology on the unitary group of the latter 
(the topology generated by the von 
Neumann algebra predual) coincides with the strong 
operator topology induced from $U_\infty$.
The result follows.
\end{proof}

Thus we have the following immediate corollary.
\begin{corollary} 
Every abelian Polish group $G$ is isomorphic with a topological factor-group
of a closed subgroup of $L^0(X,U(1))$.
\end{corollary}
\begin{proof} Since $\ell_1$ is a closed subgroup of $U_\infty$ by 
Proposition \ref{megrel},
$\ell_1$ is a closed subgroup of $L^0(X,U(1))$ by Theorem \ref{L0}.
The statement now follows from Corollary \ref{ell-on}.
\end{proof}

In fact, because we are dealing with abelian groups here, a stronger
statement ensues.

\begin{corollary} 
Every abelian Polish group $G$ is isomorphic with a closed
subgroup of a topological factor-group of $L^0(X,U(1))$.
\end{corollary}

\begin{proof} 
Let $F$ be a closed subgroup of $\ell_1$, given by Corollary \ref{ell-on},
with the property that
$G$ embeds into $l_1/F$ as a closed subgroup. 
Embed $\ell_1$ into $L^0(X,U(1))$ as a closed subgroup using
Theorem \ref{L0}. Then $G$ is isomorphic to a closed subgroup
of $L^0(X,U(1))/F$.
\end{proof}

The topological group
$L^0(X,U(1))$ seems to play an important role in 
the theory of ``large" topological
groups. Among other known properties of $L^0(X,U(1))$ are the following two.

\begin{theorem} 
$L^0(X,U(1))$ is a monothetic topological group.
\qed
\end{theorem}

(See e.g. \cite{Gl}, or else a simple proof from \cite{M-P}.) 

\begin{theorem}[\cite{Gl}; also Furstenberg and Weiss, unpublished]
The topological group
$L^0(X,U(1))$ is extremely amenable, that is,
every continuous action of $L^0(X,U(1))$ on a compact space
has a fixed point. 
\qed
\end{theorem}

Embeddability into monothetic groups is closed under taking factor-groups
(\cite{M-P}). It is also evident that a topological factor-group of an
extremely amenable group is extremely amenable. Thus combining
the previous three results, we not only
obtain a different proof of the main result of \cite{M-P}, but also 
strengthen it as follows.

\begin{theorem} 
Every separable metrizable abelian topological group embeds as
a topological subgroup into a monothetic extremely amenable 
metrizable topological group. This group is a topological 
factor-group of $L^0(X,U(1))$.
\qed
\end{theorem}

It was previously known that every topological group embeds into
an extremely amenable group \cite{Pe}, 
but an abelian version of the result
appears here for the first time.

\section{Borel actions
\label{borel}}

If a Polish group $G$ is a topological factor-group of a Polish group $H$, then
any Borel $G$-space $X$ is trivially a Borel $H$-space. 
Moreover, the orbit equivalence relations $E^X_G$ and $E^X_H$ are the same. 

The following can be found e.g. in \cite{BK}, Theorem 2.3.5.

\begin{theorem}[Mackey]
If $G$ is a closed subgroup of a Polish
group $H$, then any Borel $G$-space $X$ can be extended 
to a Borel $H$-space $Y$
such that every $G$-orbit in $X$ is contained in exactly one $H$-orbit in $Y$.
Moreover, it follows that $E^X_G\leq_B E^Y_H$. 
\qed
\end{theorem}

Combining these results, we obtain:

\begin{corollary}
If $G$ is the topological factor-group of a closed subgroup 
of $H$, then for any Borel $G$-space
$X$ there is a Borel $H$-space $Y$ such that $E^X_G\leq_B E^Y_H$.
\qed
\end{corollary}

The same conclusion holds if $G$ is a closed subgroup of 
a topological factor-group of $H$.
But for abelian Polish group $H$ these two conditions are in fact equivalent. 

Our results from previous sections thus imply Borel 
reducibility results for orbit equivalence
relations induced by actions of the groups mentioned. 
Let us first summarize the universality
results we have essentially proved.

\begin{theorem} Every abelian Polish group is 
the topological factor-group of a closed
subgroup of any of the following groups:
\begin{enumerate}
\item[(i)] the additive group of the Banach space $\ell_1$;
\item[(ii)] the additive group of the Banach space $C([0,1])$;
\item[(iii)] $L^0(X,U(1))$, where $X$ is any uncountable 
standard Borel space with a
non-atomic Borel measure;
\item[(iv)] the full unitary group $U_\infty$.
\end{enumerate}
\label{summary}
\qed
\end{theorem}

Clause (ii) follows from the well known fact 
that the Banach space $C([0,1])$ is 
a universal separable Banach space. 

One could also add on the list two other groups.
One is the projectively universal abelian Polish group 
$\widehat{A}({\mathcal N})$, where $\mathcal N$ is the Baire space of 
infinite sequences of natural numbers. 
The proof of its universality follows from 
Proposition \ref{factors} modulo known
properties of the Baire space (equipped with the 
standard ultrametric). See \cite{Ke,S-P-W}.

The other group is the injectively universal abelian Polish group 
recently constructed by Shkarin \cite{Shk}. 
Both these groups are not quite as familiar as
the universal groups we consider here.

Our main Borel reducibility result is now immediate.

\begin{theorem} 
Let $G$ be an abelian Polish group 
and $X$ be a Borel $G$-space.
Then there is a Borel $U_\infty$-space $Y$ such that 
$E^X_G\leq_B E^Y_{U_\infty}$.
\label{main}
\end{theorem}

Orbit equivalence relations induced by abelian Polish group actions are a rich
source of examples in the descriptive set theory of equivalence relations. To name
a few important equivalence relations that have been studied intensively: the
shift equivalence relations on $\R^\N$ by classical 
Banach spaces $\ell_p$ $(p\geq 1)$ or $c_0$,
the equivalence relations on ${\mathcal P}(\N)$ (the power set of $\N$) by the
natural actions of Polishable Borel ideals. 
There has been some hope that an example of an equivalence
relation not Borel reducible to any $U_\infty$-orbit equivalence relation would be among
these examples; now we know there can be none.

Of course Theorem \ref{summary} also implies that any orbit equivalence 
relation
induced by an abelian Polish group action must be Borel reducible to 
one by an action of 
either $\ell_1$, $C([0,1])$ or $L^0(X,U(1))$. 
Furthermore we can identify some universal
equivalence relations among those induced by abelian Polish group actions. 
These are
summarized in the following theorem.
\begin{theorem}
Let $G$ be either $\ell_1$, $C([0,1])$ or $L^0(X,U(1))$. 
Let ${\mathcal F}(G)$ be the
space of all closed subsets of $G$ with the Effros Borel structure. Let $G$ act on
${\mathcal F}(G)$ by multiplication and denote the orbit equivalence relation by
$E_G$. Then $E_G$ is universal among all orbit equivalence relations induced by
abelian Polish group actions, i.e., for any abelian Polish group $H$ and Borel $H$-space
$X$, $E^X_H\leq_B E_G$.
\label{action}
\end{theorem}
\begin{proof}
By Theorem 3.5.3 of \cite{BK}, it suffices to note that $G\times \Z$ is 
isomorphic to a closed subgroup of $G$.
\end{proof}

The following problems seem to be open:
Is Theorem \ref{action} true for $\ell_p$ $(p>1)$, especially $\ell_2$, and $c_0$?

\section*{Acknowledgments}
The authors are grateful to Michael Megrelishvili for most useful
and stimulating
discussions and pointers to the literature, and to
Leon Vaserstein for helpful references. 
Research of the first named author (S.G.) was
supported by United States NSF Grant DMS-0100439 and a University of North
Texas Faculty Research Grant.
Research of the second named author (V.P.) was supported by the
Marsden Fund of the Royal Society of New Zealand. 
The same author thanks the Department of 
Mathematics of the University of North Texas for hospitality extended
in November 2001 during his visit, 
which was partly supported from the Roy McLeod Millican Memorial Fund.


\begin{thebibliography}{100}

\bibitem{AE} R. Arens and J. Eells,
{\it On embedding uniform and topological spaces,}
 Pacific J. Math.{\bf 6} (1956), 397--403.

\bibitem{BK} H. Becker and A.S. Kechris,
The Descriptive Set Theory of Polish Group Actions,
London Math. Soc. Lecture Note Ser. {\bf 232}, 
Cambridge University Press, 1996.


\bibitem{CH} J.P.R. Christensen and W. Herer,
{\it On the existence of pathological submeasures and the construction of exotic
topological groups,}
Math. Ann. {\bf 213} (1975), 203--210.

\bibitem{Fl1} J. Flood,
{\it Free Topological Vector Spaces,} Ph.D. thesis, 
Australian National University, Canberra, 1975.

\bibitem{Fl2} J. Flood,
{\it Free locally convex spaces,} Dissertationes Math.
{\bf CCXXI} (1984), PWN, Warczawa.

\bibitem{Gl} S. Glasner,
{\it On minimal actions of Polish groups,} Top. Appl. {\bf 85}
(1998), 119--125.

\bibitem{G-N} H. Gl\"ockner and K.-H. Neeb,
{\it Minimally almost periodic Abelian groups and 
commutative $W^\ast$-algebras,} in: 
Mart\'\i n Peinador, E. (ed.) at al.,   
Nuclear groups and Lie groups. 
Selected lectures of the workshop, Madrid, Spain, September 1999. 
Lemgo: Heldermann Verlag. Res. Expo. Math. {\bf 24} (2001), 163--185.

\bibitem{Gr} M.I. Graev, 
{\it Theory of topological groups I,}
Uspekhi Mat. Nauk (N.S.) {\bf 5} (1950), 3--56 (in Russian).

\bibitem{Kant} L.V. Kantorovich, 
{\it On the transfer of masses,} Dokl. Akad. Nauk USSR {\bf 37} 
(1942), no. 7--8, 227--229 (in Russian).

\bibitem{K-R} L.V. Kantorovich and G.Sh. Rubinshtein,
{\it On a function space in certain extremal problems,}
Dokl. Akad. Nauk USSR {\bf 115} (1957), No. 6, 1058--1061 (in Russian).

\bibitem{Ke} A.S. Kechris,
{\it Topology and descriptive set theory,}
Top. Appl. {\bf 58} (1994), 195--222.

\bibitem{LZ} J. Lindenstrauss and L. Tzafriri,
Classical Banach Spaces. Vol. I,
Reprint of the 1977 edition, Springer--Verlag, Berlin a.o.,
1996.


\bibitem{Meg} M.G. Megrelishvili,
{\it Reflexively but not unitarily representable 
topological groups,} Topology Proceedings, to appear.
E-print available from: 
{\tt http://www.cs.biu.ac.il/$^\sim$megereli/nunit.ps}

\bibitem{mich} E. Michael,
{\it A short proof of the Arens--Eells embedding theorem,}
Proc. Amer. Math. Soc. {\bf 15} (1964), 415--416.

\bibitem{M-M} S.A. Morris and C.E. McPhail,
{\it The variety of topological groups generated by 
the class of all Banach spaces,} in:
Dikranjan, Dikran (ed.) et al., 
Abelian groups, module theory, and topology. 
Proceedings in honour of Adalberto Orsatti's 60th birthday, 
Padua, Italy, 1997. New York, NY: Marcel Dekker. 
Lect. Notes Pure Appl. Math. {\bf 201} (1998), 319--325.

\bibitem{M-P} S.A. Morris and V. Pestov,
{\it Subgroups of monothetic groups,} Journal of Group Theory
{\bf 3} (2000), 407--417.

\bibitem{N} M.A. Naimark,
Normed Algebra,  Wolters-Noordhoff Publishing, 1972.

\bibitem{Pe} V. Pestov,
 {\it Ramsey--Milman phenomenon, Urysohn metric spaces,
and extremely amenable groups.} -- Israel Journal of Mathematics,
to appear (2002). 
Available from:
{\tt http://arXiv.org/abs/math.FA/0004010}

\bibitem{RR} S.T. Rachev and L. R\"uschendorf,
Mass Transportation Problems. Volume I: Theory. Volume II:
Applications, Springer, NY--Berlin--Heidelberg, 1998.

\bibitem{Rai} D.A. Ra\u\i kov,
{\it Free locally convex spaces for uniform spaces,}
 Mat. Sb. (N.S.) {\bf 63} (1964), 582--590 (in Russian).

\bibitem{Sak} M. Sakarovitch,
{\it Linear Programming,} Springer texts in electrical engineering,
Springer-Verlag, New York, 1983. 

\bibitem{S-P-W} D. Shakhmatov, J. Pelant, and S. Watson,
{\it A universal complete metric abelian group of a given weight,}
in: Bolyai Society Math. Studies {\bf 4} (1993),
431--439 (Topology with Applications, Szeksz\'ard (Hungary), 1993).
Available from: {\tt http://at.yorku.ca/p/a/b/c/02.aim/1.htm}

\bibitem{Shk} S.A. Shkarin, 
{\it On universal abelian topological groups,} (Russian) Mat. Sb. {\bf 190} (1999), no. 7, 
127--144; translation in Sb. Math. {\bf 190} (1999), no. 7-8, 1059--1076.

\bibitem{S} I.J. Shoenberg,
{\it On certain metric spaces arising from Euclidean spaces by a change of metric and their
imbedding in Hilbert space,} 
Ann. Math. (2nd Ser.) {\bf 38} (1937), no. 4, 787--793.

\bibitem{T} M.G. Tkachenko,
{\it On completeness of free abelian topological groups,}
Soviet Math. Dokl. {\bf 27} (1983), 341--345.

\bibitem{Top} D.M. Topping, 
Lectures on Von Neumann Algebras, van Nostrand Reinhold Company, London, 1971.

\bibitem{U} V.V. Uspenski\u\i,
{\it Free topological groups of metrizable spaces,}
Math. USSR-Izvestiya {\bf 37} (1991), 657--680.

\end{thebibliography}
\end{document}